\documentclass{article}
\usepackage[dvips]{graphicx}
\usepackage{maamonthly}

\usepackage{amsmath}
\usepackage{amssymb}
\newtheorem{lemma}{Lemma}
\newtheorem{theorem}[lemma]{Theorem}
\newtheorem{puzzle}{Puzzle}
\newrmtheorem{definition}{Definition}

\newcommand{\F}{\mathbb{F}}
\newcommand{\Proj}{\mathbb{P}}
\newcommand{\meet}{\wedge}
\newcommand{\join}{\vee}
\newcommand\Gal{\mathop{\mathrm{Gal}}\nolimits}
\newcommand\Aut{\mathop{\mathrm{Aut}}\nolimits}
\newcommand\PSL{\mathop{\mathrm{PSL}}\nolimits}
\newcommand\PGL{\mathop{\mathrm{PGL}}\nolimits}
\newcommand{\sbinom}[2]{\genfrac{[}{]}{0pt}{}{#1}{#2}}

\begin{document}

\title{Projective Geometry over $\F_1$ and
the Gaussian Binomial Coefficients}
\author{Henry Cohn}
\maketitle

\markright{Projective Geometry over $\F_1$}

\section{INTRODUCTION.}

There is no field with only one element, yet there is a
well-defined notion of what projective geometry over such a field
means.  This notion is familiar to experts and plays an
interesting role behind the scenes in combinatorics and algebra,
but it is rarely discussed as such.  The purpose of this article
is to bring it to the attention of a broader audience, as the
solution to a puzzle about Gaussian binomial coefficients.

\section{GAUSSIAN BINOMIAL COEFFICIENTS.}

What form does the binomial theorem take in a noncommutative ring?
In general one can say nothing interesting, but certain special
cases work out elegantly.  One of the nicest, due to
Sch\"utzenberger \cite{Sc}, deals with variables $x$, $y$, and $q$
such that $q$ commutes with $x$ and $y$, and $yx = qxy$.\break
Then there are polynomials $\sbinom{n}{k}_q$ in $q$ with integer
coefficients such that
\begin{equation}
\label{q-binomial}
(x+y)^n = \sum_{k=0}^n \sbinom{n}{k}_q x^ky^{n-k}.
\end{equation}
These polynomials are called \textit{Gaussian\footnote{Needless to
say, Gauss discovered them in a slightly different context. See
\cite[pp.~16--17]{G1} for how they arose in his astonishing
evaluation of the quadratic Gauss sum, and \cite[p.~462]{G2} for
another version of the \mbox{$q$-binomial} theorem (this time
commutative), but keep in mind that here the dot for
multiplication has lower precedence than addition!} binomial
coefficients\/} or \textit{$q$-binomial coefficients\/}.  They can be
calculated recursively using
\begin{equation}
\label{recurrence} \sbinom{n}{k}_q = \sbinom{n-1}{k}_q +
q^{n-k}\sbinom{n-1}{k-1}_q,
\end{equation}
together with the boundary conditions $\sbinom{n}{0}_q =
\sbinom{n}{n}_q = 1$.  To see why, note that writing
\[
(x+y)^n = (x+y)^{n-1}x+(x+y)^{n-1}y
\]
and keeping careful track of how many times $y$ moves past $x$
shows that the coefficients of~(\ref{q-binomial}) satisfy the
recurrence (the boundary conditions are obvious).

Setting $q=1$ yields the ordinary binomial coefficients and
recurrence (i.e., Pascal's triangle).  The analogy between
Gaussian and ordinary binomial coefficients can be strengthened
as follows. Define the \textit{$q$-analogue\/} of the natural number
$n$ by
\[
[n]_q = 1+q+\dots+q^{n-1}
\]
(note that setting $q=1$ yields $n$) and the \textit{$q$-factorial\/}
by $[0]_q!=1$ and
\[
\vadjust{\newpage}
[n]_q! = [1]_q[2]_q\dots[n]_q
\]
for $n\ge1$. Then it is not hard to prove by induction using~(\ref{recurrence}) that
\begin{equation}
\label{q-quotient}
\sbinom{n}{k}_q = \frac{[n]_q!}{[k]_q!\,[n-k]_q!} =
\frac{[n]_q[n-1]_q\dots[n-k+1]_q}{[k]_q!},
\end{equation}
in perfect parallel with the case $q=1$.  Note that it is not at
all obvious that the right-hand side of~(\ref{q-quotient}) is a
polynomial in $q$, although that follows from the recurrence
rela-tion.

Gaussian binomial coefficients are far more than just a
construction from algebra.\break
Indeed, they arise in a startling number of combinatorial
problems. For a taste (due in this form to P\'olya \cite{Po},
although it is equivalent to a much earlier theorem on
partitions---see section~4 in \cite{Sy}), imagine an $m \times n$
box with opposite corners at\break
$(0,0)$ and $(m,n)$, where $m$ and $n$ are positive integers. It
is a standard fact of combi-\break
natorics that there are $\binom{m+n}{m}$ paths from $(0,0)$ to
$(m,n)$ made up of steps of one unit\break
up or right (each path consists of $m+n$ steps, among which one
can freely choose which $m$ go right). Let $f(m,n,a)$ be the
number of such paths that enclose area $a$\break
with the bottom and right walls of the box. Then the Gaussian
binomial coefficients are generating functions for this quantity:
\[
\sum_{a = 0}^{mn} f(m,n,a) q^a = \sbinom{m+n}{m}_q.
\]
There is a straightforward proof using~(\ref{recurrence}), but one
can also see directly how this corresponds to the $q$-binomial
theorem (a good exercise for the reader).  More details on this
interpretation and other related ones can be found in the
excellent expository ar\-ticle~\cite{PoAl}.

For our purposes, the crucial interpretation of Gaussian binomial
coefficients is given by the following theorem about linear
algebra over the finite field $\F_q$ with $q$ elements (this
theorem's early history is not fully known---see~\cite[p.~278]{Ku}
and~\cite[p.~227]{A}):

\begin{theorem}
\label{field-count}
If $q$ is a prime power, then $\sbinom{n}{k}_q$ is the number of
$k$-dimensional subspaces of\/ $\F_q^n$.
\end{theorem}

\begin{proof}
If we substitute $[n]_q = (q^n-1)/(q-1)$ into~(\ref{q-quotient}),
we find that
\begin{eqnarray*}
\sbinom{n}{k}_q &=& \frac{(q^n-1) (q^{n-1}-1) \dots
(q^{n-k+1}-1)}{(q^k-1) (q^{k-1}-1) \dots (q-1)}\\
&=& \frac{(q^n-1) (q^n-q) \dots (q^n-q^{k-1})}{(q^k-1) (q^{k}-q)
\dots (q^k-q^{k-1})}.
\end{eqnarray*}
Now consider the number of ways to choose a $k$-tuple
$(v_1,\dots,v_k)$ of linearly independent vectors in $\F_q^n$. If
we choose the vectors consecutively, then $v_1$ can be any nonzero
vector, and the only restriction on $v_i$ is that it must not be
one of the $q^{i-1}$ linear com-binations
of $v_1,\dots,v_{i-1}$. Thus, there are $q^n-q^{i-1}$ choices for
$v_i$, and
\[
(q^n-1) (q^n-q) \dots (q^n-q^{k-1})
\]
$k$-tuples total.  Each $k$-tuple spans a $k$-dimensional
subspace of $\F_q^n$, and each subspace is spanned by
\[
\vadjust{\newpage}
(q^k-1) (q^k-q) \dots (q^k-q^{k-1})
\]
$k$-tuples (a second application of the same argument, with
$n=k$). Therefore there are
\[
\frac{(q^n-1) (q^n-q) \dots (q^n-q^{k-1})}{(q^k-1) (q^{k}-q) \dots
(q^k-q^{k-1})}
\]
$k$-dimensional subspaces, as desired.\qed
\end{proof}

One can also prove Theorem~\ref{field-count} using the
recurrence~(\ref{recurrence}), but this proof is pretty. In this
form it goes back at least to \cite{GR}, and a similar proof
occurs in Burnside's 1897 group theory book \cite{Bu} (see
pages~58--60, or pages~109--111 in the second edition from 1911).

Theorem~\ref{field-count} suggests a strong analogy between
subsets of a set (the $q=1$ case) and subspaces of a vector space
(the prime power case).  This analogy extends much further, and
has been developed by numerous authors.  See, for example, the
seminal paper \cite{GR} by Goldman and Rota.  At the end of
\cite{GR}, the authors ask for an explanation of why this analogy
holds.  It's one thing to observe it in the formulas, but quite
another to describe a consistent combinatorial picture in which
subsets appear naturally as a degenerate case of subspaces.

\begin{puzzle}
\label{gaussian-puzzle} In what way is an $n$-element set like
$\F_1^n$ (and subsets like subspaces)?
\end{puzzle}

Of course, there is no field $\F_1$ with only one element, but
there is a trivial ring, and it is merely a convention that we do
not call it a field.  However, it is an excellent convention,
because the trivial ring has no nontrivial modules (if $x$ is an
element of\break
a module, then $x=1x=0x=0$).  Calling it a field would not help
solve Puzzle~\ref{gaussian-puzzle}, since $\F_1^n$ does not depend
on $n$.

I know of no direct solution to this puzzle, nor of any way to make sense
of vector spaces over $\F_1$.  Nevertheless, the puzzle can be
solved by an indirect route: it becomes much easier to understand
when it is reformulated in terms of projective geometry.  That may
not be surprising, if one keeps in mind that many topics, such as
intersection theory, become simpler when one moves to projective
geometry. (The papers \cite{Ko} and \cite{W} also shed light on
this puzzle by indirect routes, but not by using projective
geometry.)

\section{PROJECTIVE GEOMETRY.}

Recall that projective geometry is a beautiful and symmetric
completion of affine geometry.  Given any field $F$,\footnote{In
fact, any division algebra will do, but we are interested in
finite projective geometries and all finite division algebras are
fields.  This theorem was first stated by Wedderburn in \cite{MW},
but the first of his three proofs has a gap, and Dickson gave a
complete proof before Wedderburn did.  See \cite{Pa} for details.}
one can construct the $n$-dimensional \textit{projective space\/}
$\Proj^n(F)$ as the space of lines through the origin in
$F^{n+1}$.  Equivalently, points in $\Proj^n(F)$ are equivalence
classes of nonzero points in $F^{n+1}$\break
modulo multiplication by nonzero scalars. We write
$[x_0,\dots,x_n]$ for the equivalence class of $(x_0,\dots,x_n)$
(these coordinates are called \textit{homogeneous coordinates\/}).
Affine $n$-space $F^n$ is embedded into $\Proj^n(F)$ via
$(x_1,\dots,x_n) \mapsto [1,x_1,\dots,x_n]$ (these are known as
\textit{inhomogeneous coordinates\/}), and the points with
homogeneous coordinates $[0,x_1,\dots,x_n]$ form a copy of
$\Proj^{n-1}(F)$ called the set of \textit{points at infinity\/}.
Continuing this process on the points at infinity recursively
partitions $\Proj^n(F)$ into affine pieces of each dimension up to
$n$. This point of view makes projective space look asymmetric,
but of course we can see from the definition that $\Proj^n(F)$ is
completely symmetric.

Just as points in $\Proj^n(F)$ correspond to lines through the
origin in $F^{n+1}$, lines in\break
$\Proj^n(F)$ correspond to planes through the origin in $F^{n+1}$,
and in general $k$-dimensional subspaces of $\Proj^n(F)$
\vadjust{\newpage}correspond to $(k+1)$-dimensional vector
subspaces of $F^{n+1}$.  One subspace of $\Proj^n(F)$ is contained
in another if that containment holds for the corresponding vector
subspaces of $F^{n+1}$.  We identify subspaces of $\Proj^n(F)$
with the sets of points of $\Proj^n(F)$ they contain (it is easy
to check that if they contain exactly the same points, then they
are equal).  It is convenient to consider the empty set as a
$(-1)$-dimensional subspace of $\Proj^n(F)$, which is consistent
with the foregoing definition.

The points of a $k$-dimensional subspace of $\Proj^n(F)$ are
determined by $n-k$ in-dependent
linear constraints in homogeneous coordinates (the defining
equations of the corresponding vector subspace).  In terms of
inhomogeneous coordinates for the affine subspace $F^n$, these
constraints amount to $n-k$ inhomogeneous linear equations.
\mbox{Every} $k$-dimensional affine subspace of $F^n$ is the
solution set of some equations of this sort, but not all such
collections of equations have $k$-dimensional affine solution
sets: because they are inhomogeneous equations, their solution
sets in $F^n$ may have di-mension
less than $k$, or may even be empty.  In that case most points of
the projective subspace are at infinity, and its intersection with
affine space is small.

Given two subspaces $S$ and $T$ of projective space, let $S \meet
T$ (``$S$ meet $T$'') and $S \join T$ (``$S$ join $T$'') denote
their intersection and span, respectively (i.e., take the
intersection and span of the corresponding vector subspaces of
$F^{n+1}$). The meet is their greatest lower bound under
containment, and the join is their least upper bound.  Among the
most important properties of meets and joins in projective space
is the following fact of linear algebra, called the
\textit{modular law\/}:
\[
\dim(S) + \dim(T) = \dim(S \meet T) + \dim(S \join T).
\]

The modular law implies many of the familiar properties of
projective geometry.  For example, let $S$ and $T$ be two
distinct lines in $\Proj^2(F)$.  Then $S \join T = \Proj^2(F)$,
and it follows from the modular law that $\dim(S \meet T) = 0$,
(i.e., $S$ and $T$ intersect in a point).  Similarly, let $S$ and
$T$ be distinct points in $\Proj^2(F)$.  Then $\dim(S \meet T) =
-1$, and it follows that $S\join T$ is a line and thus there is a
unique line through $S$ and $T$ (unique because every subspace
containing $S$ and $T$ contains $S \join T$).

Theorem~\ref{field-count} can be trivially reformulated in terms
of projective geometry:

\begin{theorem}
\label{projective-count}
If $q$ is a prime power, then $\sbinom{n+1}{k+1}_q$ is the number
of $k$-dimensional subspaces of\/ $\Proj^n(\F_q)$.
\end{theorem}

Puzzle~\ref{gaussian-puzzle} has a projective analogue as well:

\begin{puzzle}
\label{projective-puzzle}
In what way is an $(n+1)$-element set like $\Proj^n(\F_1)$ (and
subsets like subspaces)?
\end{puzzle}

This reformulation of the puzzle is the one we will explain.  Our
goal is to make sense of projective geometry over $\F_1$. However,
it does not fit into the linear-algebraic framework in which we
have been working. Instead, we must give a more combinatorial
definition of projective geometry, which will include not only the
case $q=1$, but also some additional projective geometries we have
not yet seen.

\begin{definition}
\label{projective} A \textit{projective geometry of order $q$\/}
is a finite set $P$ (whose elements are called \textit{points\/}),
a set $L$ of subsets of $P$ (whose elements are called
\textit{subspaces\/}), and a function $\dim : L \to
\{-1,0,1,\dots\}$ satisfying the following axioms:
\begin{enumerate}
\item \label{lattice} $L$ forms a lattice when partially ordered by
containment.  In other words, each pair of elements $S$ and $T$
has a greatest lower bound $S \meet T$ and a least upper bound $S
\join T$ in $L$ under $\subseteq$.

\newpage
\item \label{increasing} The function $\dim$ is strictly increasing:
if $S$ and $T$ belong to $L$ and $S \subsetneq T$, then $\dim(S) <
\dim(T)$.

\item \label{pointsinL} For all $x$ in $P$, $\{x\}$ is a member of $L$,
as is $\emptyset$.

\item \label{calibration} For $S$ in $L$, $\dim(S) = -1$ if and only if
$S = \emptyset$, and $\dim(S)=0$ if and only if $S = \{x\}$ for
some $x$ in $P$.

\item For $S$ and $T$ in $L$,
\[
\dim(S) + \dim(T) = \dim(S \meet T) + \dim(S \join T).
\]

\item \label{line-count} If $S$ is a member of $L$ with $\dim(S)=1$,
then $|S|=q+1$.
\end{enumerate}
\end{definition}

The terminology ``of order $q$'' is unfortunate but standard.  It
does not mean that there are $q$ points; instead, think of it as
meaning that we are working over a field with $q$ elements, as in
the case of $\Proj^n(\F_q)$, although that may not be true.  We
have made no attempt to specify a minimal set of axioms.  For
example, Axiom~\ref{increasing} follows from the other axioms. It
is essentially a theorem of Birkhoff \cite{B} that these axioms
are equivalent to other standard definitions (``essentially''
because our axioms differ slightly from Birkhoff's, but the
equivalence is not hard to prove), with the exception that most
people require $q>1$ before they use the term ``projective
geometry.''

Note that it follows from Axioms~\ref{lattice} and~\ref{pointsinL}
that $P$ belongs to $L$, since the join of all the
zero-dimensional subspaces must be $P$. We define the
\textit{dimension\/} of the geometry to be $\dim(P)$.

The complete list of finite projective geometries of order greater
than one is still unknown.  Veblen and Bussey \cite{VB} used an
approach due to Hilbert \cite{H} to classify those that satisfy
the Desargues theorem (if two triangles in a plane are in
perspective from a point, then they are in perspective from a
line).  They attempted to coordinatize the geometry, and the
Desargues theorem was needed to obtain associativity;  when it
holds, the geometry must be a projective geometry over a finite
field.  The usual proof of the Desargues theorem involves lifting
to three-dimensional space, and in fact the theorem holds in every
projective geometry of dimension greater than two.  Thus, the only
finite projective geometries remaining to be classified are the
projective planes, and in particular those that cannot be
embedded into higher-dimensional spaces. Veblen and Wedderburn
\cite{VMW} constructed examples of finite projective planes that
do not satisfy the Desargues theorem and are therefore not
defined over finite fields, but a complete list is not known. All
known examples have prime power order, and only two limitations
on the order have been established: Bruck and Ryser \cite{BR}
proved if the order is $1$ or $2$ modulo $4$ then it must be a
sum of two squares, and Lam, Swiercz, and Thiel \cite{LST}
checked by a massive computer search that the order cannot be
$10$.  In particular, it is not known whether there is a
projective plane of order $12$.  It is worth pointing out that a
projective plane of order $q$ can be defined far less verbosely
than in Definition~\ref{projective}: it is a finite set of points
with certain subsets called ``lines'' such that not all the
points lie on one line, each line has $q+1$ points, each pair of
distinct points is on a unique line, and each pair of distinct
lines intersects in a unique point.  (In fact, simply requiring
that each line must have at least three points implies that they
all have the same number of points.) Classifying these objects is
a natural and important combinatorial problem.

We can now solve Puzzle~\ref{projective-puzzle} by identifying the
projective geometries of order $1$.  They are Boolean algebras:
let $L$ consist of all subsets of $P$ and set $\dim(S) = |S|-1$.
It is clear that this defines a projective geometry of order~$1$,
and it is not difficult to check using Lemmas~\ref{facts}
\vadjust{\newpage}and~\ref{point-count} that these are the only
projective geometries of order~$1$.  As desired, their subspaces
of a given dimension are counted by ordinary binomial
coefficients.

Thus, we have solved the puzzle, and seen how the Boolean algebra
of subsets of a set fits naturally as the $q=1$ case of
projective geometry.  However, to make the solution convincing, we
must give a unified proof of the analogue of
Theorem~\ref{projective-count} for every projective geometry of
order $q$.  (If the proof required case analysis, then the
apparent unification of the definition might be illusory.)
Theorem~\ref{unified} is such a unification. Before proving it, we
deduce some lemmas from the axioms of Definition~\ref{projective}.

\begin{lemma}
\label{facts} Every projective geometry $(P,L,\dim)$ of order~$q$
and dimension~$n$ has the following properties:
\begin{enumerate}
\item \label{inductive} Each element $S$ of $L$ is itself
naturally a projective geometry $(S,L',{\dim}|_{L'})$ of
order~$q$, where $L' = \{T \in L \mid  T \subseteq S\}$.

\item \label{meet} For $S$ and $T$ in $L$, $S \meet T = S \cap T$.

\item \label{distinct-unique} Every two distinct points in $P$ lie
on a unique line, and every two distinct lines intersect in at
most one point.

\item \label{build} For $S$ in $L$ and $x$ in $P$ but not in $S$,
$\dim(S \join \{x\}) = \dim(S) + 1$.

\item \label{intersect} For $S$ and $T$ in $L$ with $\dim(S) = n-1$,
either $T$ is contained in $S$ or $\dim(T \meet S) = \dim(T)-1$.
\end{enumerate}
\end{lemma}

\begin{proof}
We deal with the assertions one by one:
\begin{enumerate}
\item All of the axioms for a projective geometry hold
trivially.  Only Axiom~\ref{lattice} requires the slightest
argument: if members $T_1$ and $T_2$ of $L$ are subsets of $S$,
then $T_1 \join T_2$ is contained in $S$ by the definition of a
least upper bound, so $T_1 \join T_2$ belongs to $L'$ as desired.

\item By definition, $S \meet T \subseteq S$ and $S \meet T \subseteq T$,
so $S \meet T \subseteq S \cap T$.  On the other\break
hand, for every $x$ in $S \cap T$, $\{x\}$ is an element of $L$
that is contained in both $S$ and $T$, so $\{x\} \subseteq S \meet
T$ by the definition of the greatest lower bound.  Hence, $S \meet
T = S \cap T$.

\item This assertion follows from the modular law and
Axiom~\ref{calibration}, as in the analysis of $\Proj^2(F)$ from
earlier in the paper (except that in more than two dimensions
there can be disjoint lines).

\item We have
\begin{eqnarray*}
\dim(\emptyset) + \dim(S \join \{x\}) &=& \dim(S \meet \{x\}) +
\dim(S \join \{x\})\\
&=& \dim(S) + \dim(\{x\})\\
&=& \dim(S),
\end{eqnarray*}
from which it follows that $\dim(S \join \{x\}) = \dim(S)+1$.

\item Because $\dim(S) = n-1$, either $T$ is a subset of $S$ or
$T \join S = P$.  In the latter case,
\begin{eqnarray*}
n + \dim(T \meet S) &=& \dim(T \join S) + \dim(T \meet S)\\
&=& \dim(T)+\dim(S)\\
&=& \dim(T) + n-1,
\end{eqnarray*}
so $\dim(T \meet S) = \dim(T) - 1$.\qed
\end{enumerate}
\end{proof}

\newpage
\begin{lemma}
\label{point-count}
Every projective geometry of order $q$ and dimension $n$ contains
$[n+1]_q$ points.
\end{lemma}

\begin{proof}
We prove this by induction on $n$.  The base case $n=0$ follows
from Axiom~\ref{calibration}.  Now suppose that the lemma holds
for all dimensions less than~$n$.

By repeatedly applying assertion~\ref{build} in Lemma~\ref{facts},
one can construct a subspace $S$\break
of dimension $n-1$.  There must also be a point $x$ that does not
lie in $S$. Every line through $x$ intersects $S$ in a unique
point by assertions~\ref{meet} and~\ref{intersect}.  By
assertion~\ref{distinct-unique}, every point other than $x$ lies
on a unique line with $x$, and these lines are all disjoint except
for $x$.  Each line contains $q$ points besides $x$ by
Axiom~\ref{line-count}.  Therefore, the total number of points in
the geometry is $1+q|S| = 1+q[n]_q = [n+1]_q$, as desired ($|S| =
[n]_q$\break
by assertion~\ref{inductive} and the inductive hypothesis).\qed
\end{proof}

\begin{theorem}
\label{unified} Every projective geometry of order $q$ and
dimension $n$ contains $\sbinom{n+1}{k+1}_q$ subspaces of
dimension $k$.
\end{theorem}

\noindent (While Theorem~\ref{unified} can be proved analogously
to Theorem~\ref{field-count}, for variety we will instead use the
recurrence~(\ref{recurrence}).)

\begin{proof}
As in the preceding proof, we prove this by induction on $n$. The
base case $n=0$ is again trivial.  Thus, we suppose that the
result holds for all dimensions less than $n$.

Let $S$ be a subspace of dimension $n-1$.  By the inductive
hypothesis, there are\break
$\sbinom{n}{k+1}_q$ subspaces of dimension $k$ in $S$. By
assertion~\ref{intersect} of Lemma~\ref{facts}, every other
\mbox{$k$-dimensional} subspace intersects $S$ in a
$(k-1)$-dimensional subspace, so there are $\sbinom{n}{k}_q$
possible intersections. To complete the proof, we will show that
every $(k-1)$-dimensional subspace of $S$ extends in $q^{n-k}$
ways to a $k$-dimensional subspace not contained in $S$.

Let $T$ be a $(k-1)$-dimensional subspace of $S$.  Each extension
is of the form $T \join \{x\}$ for some $x$ not belonging to $S$
(it contains a subspace of this form and must coincide with it
because they have the same dimension), and that partitions the
complement of $S$ in $P$ into disjoint subsets, according to
whether they lie in the same extension. It follows from
Lemma~\ref{point-count} that there are $q^n$ choices of $x$
outside $S$, and that each of the extensions contains $q^k$ of
them, so there are $q^{n-k}$ extensions.  Thus, there are
\[
\sbinom{n}{k+1}_q + q^{n-k} \sbinom{n}{k}_q = \sbinom{n+1}{k+1}_q
\]
$k$-dimensional subspaces in total, as desired.\qed
\end{proof}

\section{FURTHER DIRECTIONS.}

Viewing Boolean algebra as a special case of projective geometry
can illuminate more than just the puzzle with which we started.
One interesting example, suggested by Robert Kleinberg, is the
classification of finite simple groups.  Recall that these groups
fall into four classes (see \cite[sec.~47]{Asc}).  Aside from
cyclic groups of prime order and finitely many sporadic groups,
the only finite simple groups are the simple groups of Lie type
and the alternating groups.  The simple groups of Lie type are
finite-field analogues of simple Lie groups, and it is very
reasonable to\break
expect to construct finite simple groups in this way. What may be
surprising is that the alternating groups, which to a naive
observer feel very different from the groups of Lie type, can be
brought at least partially into the same framework. In particular,
$A_n$ can be thought of as $\PSL_n(\F_1)$, as follows.

\newpage
The most basic example of a finite group of Lie type is
$\PSL_n(\F_q)$, which is simple unless $n=2$ and $q$ is $2$ or $3$
(assume from now on that we are not in these cases). It arises
geometrically as a normal subgroup of the group
$\Aut(\Proj^{n-1}(\F_q))$ of collineations of
$\Proj^{n-1}(\F_q)$ (i.e., permutations of the points that
furthermore map subspaces to subspaces).  When $n \ge 3$, the
collineation group is a semidirect product $\Gal(\F_q/\F_p)
\ltimes \PGL_n(\F_q)$ if $q$ is a power of the prime $p$ (see
Theorem~2.26 and the discussion that follows it in
\cite[pp.~88--91]{Ar}). The subgroup $\PSL_n(\F_q)$ can be derived
from $\Aut(\Proj^{n-1}(\F_q))$ by repeatedly taking the
commutator subgroup: the commutator subgroup of
$\Aut(\Proj^{n-1}(\F_q))$ is contained in $\PGL_n(\F_q)$, the
commutator subgroup of that is equal to $\PSL_n(\F_q)$, and
$\PSL_n(\F_q)$ is its own commutator subgroup because it is a
non-Abelian simple group.  (If $q$ is prime, then one needs
to take the commutator subgroup only once to reach $\PSL_n(\F_q)$.)

What should the $q=1$ analogue be?  The automorphism group of
$\Proj^{n-1}(\F_1)$ is the symmetric group $S_n$, whose
commutator subgroup is $A_n$, and $A_n$ is simple if $n \ge 5$.
This suggests that $\PSL_n(\F_1)$ should be interpreted as $A_n$.
However, it is not clear how far the analogy goes.  For example,
what happens if one sets $q=1$ in the equation
\[
|\PSL_n(\F_q)| = \frac{q^{\binom{n}{2}} (q-1)^{n-1} [n]_q!}
{\gcd(n,q-1)}
\]
(see Table~16.1 in \cite[p.~252]{Asc}, but note that $\pi$ is a
typo for $n$)?  The power of $q$ simply\break
becomes $1$, and $[n]_q!$ becomes $n!$, but the remaining factors
amount to $0/n$ rather than~$1/2$.  Is there any way to make sense
of this?  Can the analogy between $\PSL_n(\F_1)$ and $A_n$ be
extended or refined?

\begin{acknowledge}{Acknowledgments.}
I am grateful to
Robert Kleinberg, Elizabeth Wilmer, and the anonymous referees for
helpful comments on the manuscript.
\end{acknowledge}

\begin{biog}
\item[HENRY COHN] is a researcher in the theory group
at Microsoft Research.  He received his Ph.D.\ from Harvard
University in 2000 under Noam Elkies, after which he spent a year
as a postdoc at Microsoft Research before joining the
group long term in 2001.  His primary mathematical interests are
number theory, combinatorics, and the theory of computation.
\begin{affil}
Microsoft Research, One Microsoft Way, Redmond,
WA 98052-6399\\
cohn@microsoft.com
\end{affil}
\end{biog}

\vfill

\vfill\eject
\end{document}